\documentclass[aop,MSNbibl,citesort,dvips]{arximspdf}
\usepackage{eucal}
\usepackage{eufrak}
\usepackage{graphicx}

\doi{10.1214/11-AOP688}
\volume{40}
\issue{6}
\pubyear{2012}
\firstpage{2460}
\lastpage{2482}

\makeatletter

\newcommand{\mvs}{\fontfamily{mvs}\fontencoding{U}%
\fontseries{m}\fontshape{n}\selectfont}

\newcommand{\mvchr}[1]{{\mvs\char#1}}

\newcommand\Bicycle{\mvchr{174}}
\newcommand\Football{\mvchr{111}}

\newcommand{\ittexttt}{\textit} %%% buvo \texttt
\newcommand{\rmtexttt}{\mathrm} %%% buvo \texttt

\newproclaim{Remark}{Remark}
\newproclaim{Remarks}{Remarks}

\newtheorem{DCP}{Proposition}
\newtheorem{GT}[DCP]{Theorem}
\newtheorem{corollary}[DCP]{Corollary}

\newtheorem{lemma}{Lemma}

\newcommand{\implies}{\rightarrow}

\newcommand{\varlimsup}{\mathop{\overline{\lim}}}
\newcommand{\varliminf}{\mathop{\underline{\lim}}}

\newcommand{\text}[1]{\mbox{#1}}
\newcommand{\mod}{\mbox{ mod }}

\newcommand{\e}{\varepsilon}
\newcommand{\G}{\Gamma}
\newcommand{\D}{\Delta}
\newcommand{\s}{\sigma}
\newcommand{\x}{\times}
\newcommand{\C}{\mathcal C}

\newcommand{\Om}{\Omega}
\newcommand{\B}{\mathcal B}

\def\wasyfamily{\fontencoding{U}\fontfamily{wasy}\selectfont}
\def\smiley     {\mbox{\wasyfamily\char44}}
\newcommand{\Smi}{\smiley}

\newcommand{\lra}{\hbox to 1cm{\rightarrowfill}}

\makeatother

\begin{document}
\begin{frontmatter}

\title{Relative complexity of random walks in random~sceneries}
\runtitle{Relative complexity}

\begin{aug}
\author[A]{\fnms{Jon} \snm{Aaronson}\corref{}\thanksref{t1}\ead[label=e1]{aaro@post.tau.ac.il}\ead[label=u1,url]{http://www.math.tau.ac.il/\textasciitilde aaro}}
\runauthor{J. Aaronson}
\affiliation{Tel Aviv University}
\dedicated{Dedicated to the memory of Dan Rudolph}
\address[A]{School of Mathematical Sciences\\
Tel Aviv University\\
69978 Tel Aviv\\
Israel\\
\printead{e1}\\
\printead{u1}} %adresu isvedimo komanda gale!
\end{aug}

\thankstext{t1}{Supported by Israel Science Foundation Grant 1114/08.}

% HISTORY:
\received{\smonth{1} \syear{2010}}
\revised{\smonth{11} \syear{2010}}

% ABSTRACT
%
\begin{abstract}
\ittexttt{Relative complexity} measures the complexity of a probability
preserving transformation relative to a factor being a sequence of
random variables whose exponential growth rate is the relative entropy
of the extension. We prove distributional limit theorems for the
relative complexity of certain zero entropy extensions: RWRSs whose
associated random walks satisfy the $\alpha$-stable CLT ($1<\alpha\le2$).
The results give invariants for \ittexttt{relative isomorphism} of these.
\end{abstract}

% KEYWORDS
%

\end{frontmatter}

\section*{Introduction}

Invariants generalizing entropy and measuring the
``\ittexttt{com\-plexity}''
of a probability preserving transformation with zero entropy have been
introduced in~\cite{FP,F} and~\cite{KT}.

Here we consider corresponding ``\ittexttt{relative}'' notions applied to a
transformation over a factor (the classical definitions being retrieved
when the factor is trivial).

We give explicit computations of the invariants obtained
(distributional limits) for
certain \ittexttt{random walks in random sceneries}.

Let $(X,\B,m,T)$ be a probability preserving transformation, and let
\[
\frak P=\frak P(X,\B,m):=\{\text{countable, measurable partitions
of } X\}.
\]
A \textit{$T$-generator} is a partition $P\in\frak P$ satisfying $\s
(\bigcup_{n\in\mathbb Z}T^nP)=\B \mod m$.

Given the probability preserving transformation $(X,\B,m,T)$, $P\in
\frak P$ and \mbox{$n\ge1$}, the \textit{Hamming metric} on $P_n:=\bigvee
_{j=0}^{n-1}T^{-j}P$ is given by
\[
\overline d{}^{(P)}_n\bigl(a^{(1)},a^{(2)}\bigr):=\frac1n\#\bigl\{0\le k\le n-1\dvtx
a^{(1)}_k\ne
a^{(2)}_k\bigr\},
\]
where $a^{(i)}=[a^{(i)}_0,\ldots,a^{(i)}_{n-1}]=\bigcap
_{j=0}^{n-1}T^{-j}a^{(i)}_j$ $(i=1,2)$.\vspace*{1pt}

This induces the $(T,P,n)$-\textit{Hamming pseudometric} on $X$ given by
\[
d_n^{(P)}(x,y):=\overline d{}^{(P)}_n(P_n(x),P_n(y)),
\]
where $P(z)$ is defined by $z\in P(z)\in P$.

\subsection*{Relative complexity}

The following definitions are relativized versions of those in
\cite{F} and~\cite{KT}.

Given a factor $\mathcal C\subset\B$ (i.e., a $T$-invariant sub-$\s
$-algebra) and $n\ge1, \e>0$, define
$K_\C(P,n,\e)=K^{(T)}_\C(P,n,\e)\dvtx X\to\mathbb R$ by
\[
K_\C(P,n,\e)(x):=\min \biggl\{\# F\dvtx F\subset X, m\biggl(\bigcup_{z\in
F}B(n,P,z,\e)\big\|\mathcal C\biggr)(x)>1-\e\biggr\},
\]
where
\[
B(n,P,x,\e):=\bigl\{y\in X\dvtx d_n^{(P)}(x,y)\le\e\bigr\},
\]
and $m(\cdot\|\C)$ denotes conditional measure
with respect to $\C$.

Note that
\[
B(n,P,x,\e)=\bigcup_{a\in P_n\dvtx \overline
d{}^{(P)}_n(a,P_n(y))\le\e}a
\]
and is therefore a union of $P_n$-cylinders.

The random variable $K_\C(P,n,\e)$ is $\C$-measurable, and the family
\[
\{K_\C(P,n,\e)\dvtx n\ge1, \e>0\}
\]
is called the  \textit{relative complexity of $T$
with respect to $P$ given $\C$}.

The unwieldiness of this family motivates a search for one sequence
which describes its asymptotic properties. For example, one such
sequence is given as follows:

It follows from the discussions in~\cite{F} and~\cite{KT} that
{\renewcommand{\theequation}{$\bigstar$}
\begin{equation}\label{bstar}
\frac1n\log
K_\C(P,n,\e)\mathop{\lra}^{m}_{n\to\infty
, \e
\to0}
h(T,P\|\C),
\end{equation}}

\vspace*{-6pt}

\noindent where $\stackrel{m}{\lra}$ denotes convergence in measure and $h(T,P\|
\C
)$ denotes the relative entropy of the process $(P,T)$
with respect to $\C$.

We consider a distributional version amplifying this complexity
convergence in the case $h(T,P\|\C)=0$. To ``warm up'' for this we give
a sketch proof of~(\ref{bstar}) at the end of Section~\ref{sec2}.

\subsection*{Complexity sequences}

Let $(X,\B,m,T)$ be a probability preserving transformation, let $\C
\subset\B$ be a factor, let $\mathcal K=\{n_k\}_k, n_k\to\infty$
be a
subsequence and let $P\in\frak P$.\vadjust{\goodbreak}

We call the sequence $(d_k)_{k\ge1}$ $(d_k>0)$ a
\textit{$\C$-complexity sequence along $\mathcal K=\{n_k\}_k$} if
$\exists
$ a random variable $Y$ on $(0,\infty)$ such that
\renewcommand{\theequation}{a}
\begin{equation}
\frac{\log K_\C(P,n_{k},\e)}{d_{k}}
\mathop{\lra}^{{\frak d}}_{k\to\infty, \e\to0}
Y,%
\end{equation}
where $\stackrel{\frak d}{\lra}$ denotes convergence in distribution
(or just a \textit{$\C$-complexity sequence} in case $\mathcal K=\mathbb N$).

For example,
if
$h(T,P\|\C)>0$, then according to~(\ref{bstar}), $(n)_{n\ge1}$ is a
$\C
$-complexity sequence for $(T,P)$ with $Y=h(T,P\|\C)$.

We'll see (below) that if the distributional convergence (a) holds for
some $T$-generator $P\in\frak P$, then it holds $\forall$
$T$-generators $P\in\frak P$ in which case we call the sequence
$(d_k)_{k\ge1}$ a \textit{$\C$-complexity sequence for $T$ along
$\mathcal K=\{n_k\}_k, n_k\to\infty$}.

The growth rates of these are invariant under \ittexttt{relative
isomorphism} (see below).

\subsection*{Relative entropy dimension}

This is a relative, subsequence version of the \ittexttt{entropy dimension}
in~\cite{FP}.

Let $(X,\B,m,T)$ be a probability preserving transformation, let $\C
\subset\B$ be a factor and let $\mathcal K=\{n_k\}_k, n_k\to\infty$.

The \textit{upper relative entropy dimension of $T$
with respect to $\C$ along
$\mathcal K$} is
\[
\overline{\rmtexttt{E\mbox{-}dim}}_{\mathcal K}(T,\C):=\inf\biggl\{t\ge0\dvtx \frac
{\log
K(P,n_k,\e)}{n_k^t}
\mathop{\lra}^{m}_{k\to\infty, \e\to0} 0\mbox{ }\forall P\in
\frak P\biggr\}
\]
and

the \textit{lower relative entropy dimension of $T$
with respect to $\C$ along
$\mathcal K$} is
\[
\underline{\rmtexttt{E\mbox{-}dim}}_{\mathcal K}(T,\C):=\sup\biggl\{t\ge0\dvtx \exists
P\in
\frak P, \frac{\log K(P,n_k,\e)}{n_k^t}
\mathop{\lra}^{m}_{k\to\infty, \e\to0} \infty\biggr\}.
\]

In case the upper and lower entropy dimensions coincide, we call the
mutual value the \textit{relative entropy dimension of $T$
with respect to $\C$
along $\mathcal K$} and denote it by ${\rmtexttt{E\mbox{-}dim}}_{\mathcal
K}(T,\C)$.

As before, we'll drop reference to $\mathcal K$ in case $\mathcal
K=\mathbb N$ writing $\overline{\rmtexttt{E\mbox{-}dim}} (T,\C):=\overline{
\rmtexttt{E\mbox{-}dim}}_{\mathbb N}(T,\C)$ and
$\underline{\rmtexttt{E\mbox{-}dim}} (T,\C):=\underline
{\rmtexttt{E\mbox{-}dim}}_{\mathbb N}(T,\C)$.

Simple manipulation of the definitions (using the monotonicity lemma
below) shows that:

$\bullet$ if $(d_k)_{k\ge1}$ is a $\C$-complexity
sequence for $T$ along $\mathcal K=\{n_k\}_k, n_k\to\infty$,
then%\vspace*{8pt}
\renewcommand{\theequation}{\protect
\includegraphics{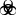}
}
\begin{equation}%[{$\hazard$}]
\label{bzard}
\overline{\rmtexttt{E\mbox{-}dim}}_{\mathcal
K}(T,\C)=\varlimsup_{k\to\infty}\frac{\log d_k}{\log n_k}
\underline{\rmtexttt{E\mbox{-}dim}}_{\mathcal
K}(T,\C)=\varliminf_{k\to\infty}\frac{\log d_k}{\log
n_k}.\hspace*{-10pt}
\end{equation}
In fact \mbox{(\ref{bzard})}
also holds
under the (more relaxed) assumption of tightness in $(0,\infty)$ of the
family $\{Z_{k,\e}:=\frac{\log K^{(T)}_\C(P,n_k,\e)}{d_k}\dvtx k\ge
1,
\e>0\}$ in the sense\vadjust{\goodbreak} that for each $\eta>0$ $\exists K_\eta\in\mathbb
N,
\e_\eta>0$ and a compact interval $J\subset(0,\infty)$ such that
\[
m([Z_{k,\e}\notin J])<\eta\qquad\forall k>K_\eta, 0<\e<\e_\eta.
\]

\subsection*{Random walk in random scenery}

A \ittexttt{random walk on random scen\-ery} (\ittexttt{RWRS}) is a skew product
probability preserving transformation, which we proceed to define in detail:

The \textit{random scenery} is an invertible, probability preserving
transformation $(Y,\C,\nu,S)$ and
the \textit{random walk on the random scen\-ery} $(Y,\C,\nu,S)$ with
\textit{jump random variable $\xi$} (assumed $\mathbb Z$-valued) is the skew product
$(Z,\B(Z),\break m,T)$
defined by
\renewcommand{\theequation}{$\heartsuit$}
\begin{equation}\label{hsuit}
Z:=\Om\x Y,\qquad m:=\mu_\xi\x\nu \quad\mbox{and}\quad
T(x,y):=(Rx,S^{x_0}y),
\end{equation}
where
\[
(\Om,\B(\Om),\mu_\xi,R):=\Bigl(\mathbb Z^\mathbb Z,\B(\mathbb Z^\mathbb Z),\prod
\operatorname{dist} \xi,\rmtexttt{shift}\Bigr)
\]
is the shift of the (independent) jump random variables.

The probability preserving transformation $(\Om,\B(\Om),\mu_\xi
,R)$ is
known as the \textit{base}.

We'll sometimes consider a corresponding RWRS with an \textit{extended
base}~$\xi$ whose base is an extension of the shift of the jumps
\[
\pi\dvtx(\Om',\B',\mu',R')\to(\Om,\B(\Om),\mu_\xi,R)
\]
and which is defined by
\renewcommand{\theequation}{$\spadesuit$}
\begin{equation}\label{ssuit}
Z':=\Om'\x Y,\qquad m':=\mu'\x\nu \quad\mbox{and}\quad
T'(x,y):=\bigl(R'x,S^{\pi(x)_0}y\bigr).
\end{equation}
The terminology RWRS was coined in~\cite{KS} where it was
attributed to Paul Shields.

There are generalizations of RWRS over more general locally compact
topological groups (not considered here) where the RWRS is constructed
using a random walk on such a group and whose scenery is a probability
preserving action of the group; see~\cite{dHS,Ba}.

As shown in~\cite{Me}, a RWRS is a K-automorphism if the random walk is aperiodic
and the scenery is ergodic.

If the scenery has finite entropy and the random walk is recurrent,
then the RWRS has the same entropy as its base.

Possibly the best known RWRS is Kalikow's $[T-T^{-1}]$ transformation,
shown in~\cite{Ka} to be not Bernoulli. For a review of this and subsequent
work on the Bernoulli properties of RWRSs, see~\cite{dHS}.

The \ittexttt{one-sided} RWRS (\ittexttt{defined} \ittexttt{as}
\ittexttt{above} \ittexttt{but} \ittexttt{with}
$\Om$ \ittexttt{replaced} \ittexttt{by} \ittexttt{the}
\ittexttt{one-sided} \ittexttt{shift} $\Om_+=\mathbb Z^\mathbb N$) is
considered, for example, in~\cite{HHR} and~\cite{Ba} where invariants for
isomorphism and the induced \ittexttt{cofiltrations} are studied.

A random walk is called $\alpha$-\textit{stable} ($\alpha\in(0,2]$) if its jump
random variable is $\alpha$-\textit{stable} in the sense that for some
\textit{normalizing\vadjust{\goodbreak}
constants} $a(n)$ (necessarily $\frac1\alpha$-regularly varying)
\[
\frac{S_n}{a(n)}\stackrel{\frak d}\lra Y_\alpha,
\]
where\vspace*{1pt} $Y_\alpha$ has the \ittexttt{standard, symmetric $\alpha$-stable} (S$\alpha$S)
distribution of order $\alpha$ on $\mathbb R$ [defined by $\mathbb E(e^{itY_\alpha
})=e^{-{t^\alpha}/\alpha}$].
A RWRS is called $\alpha$-\textit{stable} if its corresponding random walk is
$\alpha$-stable.

We see that for an extended base RWRS $T$ whose corresponding random
walk is aperiodic and $\alpha$-\textit{stable} ($\alpha\in(1,2]$):

$\bullet$ the normalizing constants $a(n)$ form a \ittexttt{Base}-complexity
sequence for $T$;

$\bullet$ $\rmtexttt{E\mbox{-}dim} (T,\rmtexttt{base})=\frac1\alpha$.

\subsection*{Organization of the paper}

We state the results more precisely in Section~\ref{sec1}. The results on
abstract relative complexity are proved in Section~\ref{sec2}. In
Section~\ref{sec3}, we
collect some random walk convergence results necessary for the proof of
the distributional convergence of relative complexity for RWRS which is
done in Section~\ref{sec4}.

%s1 ###
\section{Results}\label{sec1}
\begin{DCP}[(Distributional compactness
proposition)]\label{pi1}
For any $P\in\frak P$, $d_k>0, n_k\to\infty, \exists k_\ell\to
\infty
$ and a random variable $Y$ on $[0,\infty]$ such that
\renewcommand{\theequation}{a}
\begin{equation}
\frac{\log K_\C(P,n_{k_\ell},\e)}{d_{k_\ell}}
\mathop{\lra}^{{\frak d}}_{\ell\to\infty, \e\to0}
Y.
\end{equation}
\end{DCP}
\begin{GT}[(Generator theorem)]\label{pi2}
\textup{(a)}  If there is a countable $T$-generator $P\in\frak P$ satisfying
\[
\frac{\log K_\C(P,n_k,\e)}{d_k}
\mathop{\lra}^{{\frak d}}_{k\to\infty, \e\to0}
Y,
\]
where $Y$ is a random variable on $[0,\infty]$, then
\renewcommand{\theequation}{\protect
\includegraphics{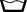}
}
\begin{equation}\label{WCotton}
\frac{\log K_\C(Q,n_k,\e)}{d_k}
\mathop{\lra}^{{\frak d}}_{k\to\infty, \e\to0}
Y\qquad \forall T\mbox{-}\rmtexttt{generators}\  Q\in\frak P.
\end{equation}

\textup{(b)} If
\[
\frac{\log K_\C(P,n_k,\e)}{n_k^t}
\mathop{\lra}^{m}_{n\to\infty, \e\to0} 0
\]
for some $T$-generator $P\in\frak P$, then
$\overline{\rmtexttt{E\mbox{-}dim}}
_{\mathcal K}(T,\C)\le t$.

\textup{(c)} If
\[
\frac{\log K_\C(P,n_k,\e)}{n_k^t}
\mathop{\lra}^{m}_{n\to\infty, \e\to0}
\infty
\]
for some $P\in\frak P$, then $\underline{\rmtexttt{E\mbox{-}dim}}_{\mathcal
K}(T,\C
)\ge t$.\vadjust{\goodbreak}
\end{GT}

We'll abuse notation by abbreviating \mbox{(\ref{WCotton})} by
\renewcommand{\theequation}{\protect
\includegraphics{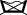}
}
\begin{equation}\label{Dwash}
\frac1{d_k}\log K^{(T)}_\C(n_k)\approx Y
\end{equation}
as in, for example,
\[
\frac1n\log K^{(T)}_\C(n)\approx h(T\|\C).
\]
\begin{GT}[(Distributional convergence theorem)]\label{pi3}
Let $(Z,\B(Z),m,T)$ be an extended base RWRS with $\alpha$-stable,
aperiodic jumps ($\alpha>1$)
and ergodic scenery $(Y,\C,\nu,S)$ satisfying $0<h(S)<\infty$, then
%
%e1
\renewcommand{\theequation}{\arabic{equation}}
\setcounter{equation}{0}
\begin{equation}
\frac1{a(n)}\log K^{(T)}_{\B(\Om)\x Y}(n)
\approx
\rmtexttt{Leb} (B_\alpha([0,1]))\cdot h,
\end{equation}
where $h:=h(S)$, $(a(n))_{n\ge1}$ are the normalizing constants of the
random walk, $\rmtexttt{Leb}$ denotes Lebesgue measure on $\mathbb R$ and
$B_\alpha$ is the S$\alpha$S process (see below).
\end{GT}

Thus, as advertised at the end of Section 0, $(a(n))_{n\ge1}$ is a
\ittexttt{base}-complexity sequence for $T$ and $\rmtexttt{E\mbox{-}dim} (T,
\rmtexttt{base})=\frac1\alpha$.

\subsection*{Relative isomorphism over a factor}

We say that the probability preserving transformations $(X_i,\B
_i,m_i,T_i)$ ($i=1,2$) are \textit{relatively isomorphic over the factors
$\C_i\subset\B_i$  $(i=1,2)$} if there is an isomorphism $\pi
\dvtx(X_1,\B
_1,m_1,T_1)\to(X_2,\B_2,m_2,T_2)$ satisfying $\pi\C_1=\C_2$.
\begin{corollary}[(Relative isomorphism corollary)]\label{pi4}
If the probability preserving transformations $(X_i,\B
_i,m_i,T_i)$ $(i=1,2)$ are relatively isomorphic over the factors $\C
_i\subset\B_i$ $(i=1,2)$, then $\forall \mathcal K=\{n_k\dvtx k\ge1\}
, d_k>0$,
\renewcommand{\theequation}{\arabic{equation}}
\setcounter{equation}{1}
\begin{equation}
\frac1{d_k}\log K^{(T_1)}_{\C_1}(n_k)\approx Y
\quad\Leftrightarrow\quad
\frac1{d_k}\log K^{(T_2)}_{\C_2}(n_k)\approx Y;
\end{equation}
\begin{eqnarray}
\overline{\rmtexttt{E\mbox{-}dim}}_{\mathcal
K}(T_1,\C
_1)&=&\overline{\rmtexttt{E\mbox{-}dim}}_{\mathcal
K}(T_2,\C_2),\nonumber\\[-8pt]\\[-8pt]
\underline{\rmtexttt{E\mbox{-}dim}}_{\mathcal K}(T_1,\C_1)&=&
\underline{\rmtexttt{E\mbox{-}dim}}_{\mathcal
K}(T_2,\C_2).\nonumber
\end{eqnarray}
\end{corollary}
\begin{corollary}[(Relative isomorphism of RWRSs)]\label{pi5}
Suppose that the aperiodic, stable, extended base RWRSs $(Z_i,\B_i,m_i,T_i)$
$(i=1,2)$ have sceneries with positive finite entropy and are relatively
isomorphic over their bases.

If $(Z_1,\B_1,m_1,T_1)$ has $\alpha$-stable jumps, then so does $(Z_2,\B
_2,m_2,T_2)$ and
\[
a^{(2)}(n)h\bigl(S^{(2)}\bigr)
\mathop{\sim}_{n\to\infty}
a^{(1)}(n)h\bigl(S^{(1)}\bigr),
\]
where $a^{(i)}$ denotes the sequence of normalizing constants of the random
walk associated to $T_i$ $(i = 1, 2)$.\vadjust{\goodbreak}
\end{corollary}
%
%Here (and throughout) $a_n\sim b_n$ means $\frac{a_n}{b_n}
%{\,\displaystyle \mathop{\lra}_{n\to\infty}\,} 1$.

%s2 ###
\section{Relative complexity}\label{sec2}
In this section, we prove Proposition~\ref{pi1}, Theorem~\ref{pi2} and
Corollary~\ref{pi4} which are relative versions
of results appearing in~\cite{KT,F} and~\cite{FP} (see the remark after the
proof of Proposition~\ref{pi1}).

\begin{pf*}{Proof of Proposition \protect\ref{pi1}
\textup{(the
distributional compactness
proposition)}}
Define $F_k\dvtx \mathbb R_+\x(0,1)\to[0,1]$ by $F_k(q,\e):=E(\exp[\frac
{-q\log K_\C(P,n_k,\e)}{d_k}])$, then $F_k(q,\e)\le F_k(q',\e')$
whenever $q\ge q', \e\le\e'$.

By Helly's theorem and diagonalization, $\exists$:

$\bullet$ a countable set $\G\subset(0,1)$;

$\bullet$ $F\dvtx\mathbb Q_+\x(0,1)\to[0,1]$ such that $F(q,\e)\le F(q',\e')$ whenever
$q\ge q', \e\le\e'$;
and a subsequence $k_\ell\to\infty$ such that
\[
F_{k_\ell}(q,\e)
\mathop{\lra}_{\ell\to\infty} F(q,\e) \qquad\forall \e\in
(0,1)\setminus\G
, q\in\mathbb Q_+.
\]

By the monotonicity of $F, F(q,\e)\downarrow F(q)$ as ${\e
\downarrow
0}$, whence
\[
F_{k_\ell}(q,\e)
\mathop{\lra}_{\ell\to\infty, \e\to0}
F(q) \qquad\forall  q\in\mathbb Q_+.
\]
Thus $\exists$ a random variable $Y$ on $[0,\infty]$ such that $F(q)=E(e^{-qY})$
and
\[
\frac{\log K_\C(P,n_{k_\ell},\e)}{d_{k_\ell}}
\mathop{\lra}^{{\frak d}}_{\ell\to\infty, \e\to0}
Y.
\]
\upqed\end{pf*}
\begin{Remark*}
To see a connection with definition 2 in~\cite{F}, note that it follows from
Proposition~\ref{pi1} (in the deterministic case) that for $B_n>0$,
the set
\[
\frak L:=\biggl\{C\in[0,\infty]\dvtx\exists n_k\to\infty, \frac{\log
K(P,n_k,\e)}{B_{n_k}}
\mathop{\lra}_{k\to\infty, \e\to0+}
C\biggr\}\ne\varnothing
\]
and
\[
\lim_{\e\to0+}\varlimsup_{n\to\infty}\frac{\log K(P,n,\e
)}{B_n}=\sup
\frak L.
\]
We turn next to the proof of the generator Theorem~\ref{pi2}.
\end{Remark*}
\begin{Monotonicitylemma*}
Suppose that $P, Q\in\frak P$ are
countable partitions such that $P\prec Q$, and suppose that
\[
\frac{\log K_\C(P,n_k,\e)}{d_k}
\mathop{\lra}^{{\frak d}}_{k\to\infty, \e\to0}
Y \quad\mbox{and}\quad  \frac{\log K_\C(Q,n_k,\e)}{d_k}
\mathop{\lra}^{{\frak d}}_{k\to\infty, \e\to0}Z;
\]
then
$Y\le Z$ [\ittexttt{in the sense that} $E(e^{-tY})\ge
E(e^{-tZ})$
$\forall t>0$].
\end{Monotonicitylemma*}
\begin{pf}
$P\prec Q \implies K_\C(P,n,\e)\le K_\C(Q,n,\e)$.
\end{pf}
\begin{lemma}\label{lemma1}
For $k\ge1, \e>0, x\in X$ and large $n\ge1$,
\[
K_\C(P_k,n,2\e)(x)\le K_\C\biggl(P,n,\frac{\e}k\biggr)(x)\le
K_\C
\biggl(P_k,n,\frac{\e}2\biggr)(x),
\]
where $P_k:=\bigvee_{j=0}^{k-1}T^{-j}P$.
\end{lemma}
\begin{pf}
Calculation shows that
$d_n^{(P_k)}(x,y)=kd_n^{(P)}(x,y)\pm\frac{k^2}n$
whence for large $n$,
$B(n,P_k,x,\frac{\e}2)\subseteq B(n,P,x,\frac{\e}k)\subseteq
B(n,P,x,2\e)$.
\end{pf}
\begin{lemma}\label{lemma2}
Let $P=\{P_n\}_{n\ge1}, Q=\{Q_n\}_{n\ge1}\in\frak P$ be ordered
partitions with
$\sum_{n\ge1}m(P_n\D Q_n)<\delta$, then $\forall \e>0, \exists
N $ such that $
\forall n\ge N$,
\[
m\bigl(\{x\in X\dvtx K_\C(Q,n,\e)(x)\ge K_\C(P,n,2\e
+2\delta)(x)\}
\bigr)>1-\e.
\]
\end{lemma}
\begin{pf}
Define $N_P, N_Q\dvtx X\to\mathbb N$ by $x\in P_{N_P(x)}, x\in
Q_{N_Q(x)}$.

By the ergodic theorem, for a.e. $x\in X, n\ge1$ large
\[
\frac1n\#\{0\le k\le n-1\dvtx N_P(T^kx)\ne N_Q(T^kx)\} =\frac1n\sum
_{k=0}^{n-1}1_\D(T^kx)<\delta,
\]
where $\D:=\bigcup_{n\ge1}P_n\D Q_n$.
It follows that $\forall \e>0, \exists N\ge1$ and sets $A_n\in
\B$ $(n\ge N)$ so that for $n\ge N$:

$\bullet$ $m(A_n)>1-\e$;

$\bullet$ $m(A_n\|\C)(x)>1-\e$ $\forall x\in A_n$;

$\bullet$ $d_n^{(Q)}(x,y)< d_n^{(P)}(x,y)+2\delta$ $\forall x, y\in A_n$.

Thus
\[
B(n,Q,x,r)\cap A_n\subseteq B(n,P,x,r+2\delta)
\qquad\forall
x\in A_n, r>0.
\]
Now fix $x\in A_n$, and suppose that $F\subset X, |F|=K_\C(Q,n,\e)(x)$
and $m(\bigcup_{z\in F}B(n$, $Q,z,\e)\|\C)(x)>1-\e$.

Let $F_1:=\{z\in F\dvtx B(n,Q,z,\e)\cap A_n\ne\varnothing\}$, and for
$z\in F_1$, choose $z'\in B(n,Q,z,\e)\cap A_n$, then
\[
\bigcup_{z\in F_1}B(n,Q,z',2\e)\supset\bigcup_{z\in F}B(n,Q,z,\e
)\setminus A_n^c.
\]

On the other hand,
\[
\bigcup_{z\in F_1}B(n,Q,z,2\e)\cap A_n\subset\bigcup_{z\in
F_1}B(n,P,z,2\e+2\delta),
\]
whence for $x\in A_n$,
\[
m\biggl(\bigcup_{z\in F_1}B(n,P,z,2\e+2\delta)\big\|\C\biggr)(x)>1-2\e
\]
and
\begin{eqnarray*}K_\C(Q,n,\e)(x)\ge|F_1|\ge K_\C(P,n,2\e+2\delta
)(x).
\end{eqnarray*}
\upqed\end{pf}

\begin{pf*}{Proof of Theorem \protect\ref{pi2} \textup{(Generator theorem)}}
We only prove (a), the proofs of (b) and (c) being analogous.

To prove\vspace*{1pt} (a), we show that every subsequence of $\{n_k\}$
has a sub-subsequence (also denoted $\{n_k\}$) along which $\frac{\log K_\C(Q,n_k,\e)}{d_k}
{\,\displaystyle \mathop{\lra}^{{\frak d}}_{k\to\infty, \e\to0}\,}Y$.

Fix a subsequence.
By Proposition~\ref{pi1}, $\exists$ a random variable $Z$ on $[0,\infty]$ and
a sub-subsequence
along which
\[
\frac{\log K_\C(P,n_k,\e)}{d_k}
\mathop{\lra}^{{\frak d}}_{k\to\infty, \e\to0}Y \quad\mbox{and}\quad \frac
{\log K_\C
(Q,n_k,\e)}{d_k}
\mathop{\lra}^{{\frak d}}_{k\to\infty, \e\to0}
Z.
\]
It suffices to show that $E(e^{-tY})= E(e^{-tZ})$ $\forall t>0$. We'll
show that $E(e^{-tY})\le E(e^{-tZ})$ $\forall t>0$ (the reverse
inequality following by symmetry).

To this end, fix $t>0, \e>0$.

$\bullet$ First choose $\kappa_0\ge1$ and $\delta>0$ such that $\forall k\ge
\kappa
_0, 0<r<\delta$
\[
E\biggl(\exp\biggl[\frac{-t\log K_\C(P,n_k,r)}{d_k}\biggr]\biggr)=E(e^{-tY})\pm\e
\]
and
\[
E\biggl(\exp\biggl[\frac{-t\log K_\C(Q,n_k,r)}{d_k}\biggr]\biggr)=E(e^{-tZ})\pm\e.
\]

$\bullet$ Next, for $0<r<\delta, \exists N=N_r\ge1, Q^{(r)}\prec P_N$ with
\[
\sum_{j\ge1}m\bigl(Q^{(r)}_j\D Q_j\bigr)<r.
\]

$\bullet$ Using Lemma~\ref{lemma2} $\exists \kappa_r>\kappa_0$ such that
\begin{eqnarray*}
&&E\biggl(\exp\biggl[\frac{-t\log K_\C(Q^{(r)},n_k,r)}{d_k}\biggr]\biggr)
\\
&&\qquad<E\biggl(\exp\biggl[\frac
{-t\log
K_\C(Q,n_k,4r)}{d_k}\biggr]\biggr)+\e\qquad\forall k\ge\kappa_r.
\end{eqnarray*}
Using Lemma~\ref{lemma1}, $\exists K_r>\kappa_r$ such that for $k>K_r$,
\begin{eqnarray*}
E\biggl(\exp\biggl[\frac{-t\log K_\C(Q^{(r)},n_k,r)}{d_k}\biggr]\biggr)
&\ge&
E\biggl(\exp\biggl[\frac{-t\log K_\C(P_N,n_k,r)}{d_k}\biggr]\biggr)\\
&\ge& E\biggl(\exp\biggl[\frac
{-t\log K_\C(P,n_k,{r}/({2N}))}{d_k}\biggr]\biggr)\\
&=&E(e^{-tY})-\e.
\end{eqnarray*}
Thus
$E(e^{-tY})\le E(e^{-tZ})+3\e$ $\forall \e, t>0$.

As mentioned above, this proves Theorem~\ref{pi2}(a).\vadjust{\goodbreak}
\end{pf*}

We note that Corollary~\ref{pi4} follows immediately from Theorem~\ref{pi2}.
\begin{pf*}{Proof sketch of (\protect\ref{bstar})}
Set
\begin{eqnarray*}
\Pi_n(x)&:=&\{a\in P_0^{n-1}(T)\dvtx m(a\|\mathcal C)(x)>0\},
\\
\Phi_{n,\e}(x)&:=&\min \biggl\{|F|\dvtx F\subset\Pi_n(x), m\biggl(\bigcup_{a\in
F}a\|
\mathcal C\biggr)(x)>1-\e\biggr\},
\\
\mathcal Q(P,n,\e)(x)&:=&\max \bigl\{\#\{c\in\Pi_n(x)\dvtx \overline d_n(a,c)\le\e
\}\dvtx
a\in\Pi_n(x)\bigr\},
\end{eqnarray*}
where
$\overline d_n$ is $(T,P,n)$-Hamming distance on $P_n$,
then
\renewcommand{\theequation}{\Bicycle}
\begin{equation}\label{bicycle}
\frac{\Phi_{n,\e}(x)}{\mathcal Q(P,n,\e)(x)}\le
K^{(T)}_{\mathcal C}(P,n,\e)(x)\le\Phi_{n,\e}(x).
\end{equation}
By the Shannon--MacMillan--Breiman theorem~\cite{Bre}, a.s.,
as $n\to\infty$
\[
I(P_n\|\C)(x)=\log\frac1{m(P_n(x)\|\C)(x)}=h(T,P\|\C)n\bigl(1+o(1)\bigr),
\]
whence by a standard counting argument, a.s.,
\renewcommand{\theequation}{\Football}
\begin{equation}\label{football}
\frac1{n}\log_2\Phi_{n,\e}
\mathop{\lra}_{n\to\infty, \e\to0} h(T\|\mathcal C).
\end{equation}
By direct estimation,
\renewcommand{\theequation}{\protect
\includegraphics{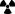}
}
\begin{equation}\label{radioactivity}
\frac1n\log\mathcal Q(P,n,\e)(x)\le\frac
1n\log
\left(|P|^{\e n}
\pmatrix{n\cr\e n}\right)
\mathop{\lra}_{n\to\infty, \e\to0} 0,
\end{equation}
whence by~(\ref{bicycle}),
\renewcommand{\theequation}{$\bigstar$}
\begin{equation}
\frac1n\log K_\C(P,n,\e)
\mathop{\lra}^{m}_{n\to\infty, \e\to0}h(T,P\|\C).
\end{equation}
\upqed\end{pf*}

The proof of Theorem~\ref{pi3} in Section~\ref{sec4} is also via~(\ref{bicycle})
and versions of
(\ref{football}) and \mbox{(\ref{radioactivity})} where the denominators $n$ are
replaced by the a sequence of normalizing constants of the random walk
$a(n)=o(n)$.

More information on random walk is needed and developed in Section~\ref{sec3}.
The version of~(\ref{football}) is established (essentially as in Section 7
of~\cite{AP}) using Skorohod's invariance principle and properties of the
range of the random walk. The proof of the \mbox{(\ref{radioactivity})} analogue uses
the invariance principle for local time as well.

%s3 ###
\section{Random walks}\label{sec3}
In this section we consider the random walk limit theorems we need to
prove Theorem~\ref{pi3}.

These are consequences of the \ittexttt{weak invariance principle} and an
\ittexttt{invariance principle for local time} as in Borodin's theorem
(below); and the properties of limit processes involved.\vadjust{\goodbreak}

$\bullet$ As in~\cite{GS}, the \textit{S$\alpha$S process} $B_\alpha$ (for $0<\alpha\le2$) is
a random function in $D([0,1])$, the \ittexttt{Donsker space} of
\textsc{cadlag} functions (Polish when equipped with the Skorokhod metric, see
\cite{Bi}) with independent, S$\alpha$S distributed increments ($B_2$~is aka
\ittexttt{Brownian motion}).

$\bullet$ The \ittexttt{weak invariance principle} says that for a $\alpha$-stable,
random walk $S_n=\sum_{k=1}^n\xi_k$,
\[
B_{\xi,n}\stackrel{\frak d}\lra B_\alpha
\]
in $D([0,1])$
where
\[
B_{\xi,n}(t):=\frac1{a_{\xi}(n)}S_{[nt]},
\]
and $a_{\xi}(n)$ are the normalizing constants (of the random walk)
satisfying
\[
\frac{S_n}{a_{\xi}(n)}\stackrel{\frak d}\lra Y_\alpha,
\]
where $\mathbb E(e^{itY_\alpha})=e^{-{|t^\alpha|}/\alpha}$.

See~\cite{D} for the case $\alpha=2$ and~\cite{GS} for $0<\alpha\le2$.

\subsection*{Local time}

For $1<\alpha\le2$, the \textit{local time at $x\in\mathbb R$} of the S$\alpha$S
process $B_\alpha$ is defined by
\[
L_\alpha(t,x):=\lim_{\e\to0+}\frac1{2\e}\int_0^t1_{[x-\e,x+\e
]}(B_\alpha(s))\,ds
\]
the limit being known to exists a.s. As shown in~\cite{Boy},
a.s., $ L_\alpha\in C_0([0,1]\x\mathbb R)$, the space of continuous functions
on $[0,1]\x\mathbb R$ tending to zero at infinity, which is Polish when
equipped with the $\sup$-norm.

We need more information about the unit range $B_\alpha([0,1])$ of the
S$\alpha
$S process.
\begin{lemma}[{\cite{EK}}]\label{lemma3}
With probability $1$, $B_\alpha([0,1])$ is Riemann integrable in $\mathbb R$
and $ L_\alpha(1,x)>0$ for $\rmtexttt{Leb}$-a.e. $x\in B_\alpha
([0,1])$.
\end{lemma}
\begin{Remark*}
More is true when $\alpha=2$. Brownian motion
$B_2$ is a.s. continuous whence $B_2([0,1])=[\min_{t\in
[0,1]}B_2(t),\max_{t\in[0,1]}B_2(t)]$. The Ray--Knight theorem
\cite{Kn,Ra}, states that a.s.,
$L_2(1,x)>0$ iff $x\in B_2([0,1])^o=(\min_{t\in[0,1]}B_2(t)$, $\max
_{t\in[0,1]}B_2(t))$.

Here and throughout, we denote the \textit{interior} (maximal open subset)
of $F\subset\mathbb R$ by $F^o$.

It is an interesting question as to whether this version of the
Ray--Knight theorem persists for $1<\alpha<2$, that is, whether
$B_\alpha([0,1])^o=\{x\in\mathbb R\dvtx L_\alpha(1,x)>0\}$ with probability
$1$.
\end{Remark*}
\begin{pf*}{Proof of Lemma~\ref{lemma3}}
By continuity of $x\mapsto L_\alpha(1,x)$, a.s.
\[
B_\alpha([0,1])^o\supset\{x\in\mathbb R\dvtx L_\alpha(1,x)>0\},\vadjust{\goodbreak}
\]
and it suffices to show that with probability $1$,
\[
\rmtexttt{Leb} \bigl(B_\alpha([0,1])\cap\{x\in\mathbb R\dvtx L_\alpha(1,x)=0\}\bigr)=0.
\]

To see this, for $F\in D([0,1])$ and $y\in[0,1]$, define
\begin{eqnarray*}L(F)(y):\!&=&
\varlimsup_{n\to\infty}2n\int_0^11_{(y-1/n,y+
1/n)}(F(t))\,dt\\ &=&
\varlimsup_{n\to\infty}2n\int_0^11_{(F(t)-1/n,F(t)+1/n)}(y)\,dt
\end{eqnarray*}
and define for $F\in D([0,1]), t\in[0,1]$
\[
\Phi(F,t):=(F,F(t),L(F)(F(t)))\in D([0,1])\x\mathbb R^2.
\]
We claim that $\Phi\dvtx D([0,1])\x[0,1]\to D([0,1])\x\mathbb R^2$ is Borel
measurable.

To see this, note first that $F\dvtx[0,1]\to\mathbb R$ is bounded, Borel
measurable ($F\in D([0,1])$ being a uniform limit of step functions), whence
$L(F)\dvtx\mathbb R\to\mathbb R$ is bounded, Borel measurable. Thus $\Phi$ is
Borel measurable.

Next, we claim that
\[
\mathbf{A}:=\{(F,y)\in D([0,1])\x\mathbb R\dvtx y\in F([0,1]),  L(F)(y)=0\}
\]
is an analytic set in $D([0,1])\x\mathbb R$.

This is because
\[
\mathbf{A}=\Pi\bigl(\Phi\bigl(D([0,1])\x[0,1]\bigr)\cap\bigl(D([0,1])\x\mathbb
R\x\{0\}\bigr)\bigr),
\]
where $\Pi\dvtx D([0,1])\x\mathbb R^2\to D([0,1])\x\mathbb R$ is the projection
$\Pi(F,x,y):=(F,x)$.

Thus $\mathbf{A}$ is $\rmtexttt{Prob}\x\rmtexttt{Leb}$-Lebesgue
measurable in $D([0,1])\x\mathbb R$
where $\rmtexttt{Prob}:=\rmtexttt{dist}\, B_\alpha\in\mathcal P(D([0,1]))$.

Next, a.s., $L(B_\alpha)(y)=L_\alpha(1,y)$ and (see, e.g.,~\cite{Ke} and references
therein) \mbox{$\forall y\in\mathbb R$},
\[
\rmtexttt{Prob}\bigl(\bigl[B_\alpha([0,1])\ni y \mbox{ and }  L_\alpha(1,y)=0\bigr]\bigr)=0.
\]
Thus, using Fubini's theorem,
\begin{eqnarray*}
&&
E\bigl(\rmtexttt{Leb} \bigl(B_\alpha([0,1])\cap[ L_\alpha(1,\cdot
)=0]\bigr)\bigr) \\
&&\qquad=\rmtexttt{Prob}\x\rmtexttt{Leb}(\mathbf{A})\\
&&\qquad=\int_{\mathbb R}E\bigl(1_{B_\alpha([0,1])}(y)1_{[ L_\alpha(1,y)=0]}\bigr)\,dy
\\
&&\qquad=
\int_{\mathbb R}\,\rmtexttt{Prob}\bigl(\bigl[B_\alpha([0,1])\ni y \mbox{ and }  L_\alpha
(1,y)=0\bigr]\bigr)\,dy\\
&&\qquad=0.
\end{eqnarray*}
\upqed\end{pf*}

\subsection*{Random walk local time}

The \textit{local time of the random walk} is
\[
N_{n,k}(x):=\#\{0\le j\le n-1\dvtx S_j(x)=k\}\qquad(n\ge1, k\in\mathbb
Z, x\in\Om).\vadjust{\goodbreak}
\]

We define the linear interpolation of $N$ by
\begin{eqnarray*}
\widehat N(n+s,k+t)&:=&
(1-s)(1-t)N_{n,k}+s(1-t)N_{n+1,k}\\
&&{}+(1-s)tN_{n,k+1}+stN_{n+1,k+1}
\end{eqnarray*}
for $s,t\in[0,1], n\in\mathbb N, k\in\mathbb Z$,
and let
\[
L_{\xi,n}(t,x):=\frac1{\overline a_{\xi}(n)}\widehat N(nt,a_{\xi}(n)x),
\]
where $\overline a_{\xi}(x):=\int_0^x\frac1{a_{\xi}(t)}\wedge1\,dt$.
\begin{Remarks*}
(i) Since $a_{\xi}(x)$ is $\frac1\alpha$-regularly varying, we have
$\overline
a_{\xi}(x)\sim\frac{\alpha}{\alpha-1}\frac{x}{a_{\xi}(x)}$.

(ii) $ L_{\xi,n}\in C_0([0,1]\x\mathbb R)$.
\end{Remarks*}

\subsection*{\texorpdfstring{Borodin's theorem~\cite{Bo1,Bo2}}{Borodin's theorem [4, 5]}}

Suppose that $(S_1,S_2,\ldots)$ is an aperiodic, $\alpha$-stable random
walk on $\mathbb Z$ with $1<\alpha\le2$, then
\[
(B_{\xi,n}, L_{\xi,n})\mathop{\lra}^{{\frak d}}_{n\to\infty}
(B_\alpha, L_\alpha
) \qquad\text{in }  D([0,1])\x C_0([0,1]\x\mathbb R).
\]

Borodin's theorem strengthens the invariance principle for local time
in Section~2 of~\cite{KS}.

Next we state the main lemma of this section. To this end we first
establish some notation.

\section*{Convergence in distribution via convergence in measure}

For the rest of this section, we'll fix $(S_1,S_2,\ldots)$, an
aperiodic, $\alpha$-stable random walk on $\mathbb Z$ with $1<\alpha\le2$ defined
on $(\Om,\B(\Om),\mu)$ as before and use the following (seemingly
stronger but) equivalent ``coupling version'' of Borodin's theorem.

Let
$(\bolds{\Om},\bolds{\mathcal{F}}):=\Om\x(D([0,1])\x C_0([0,1]\x
\mathbb
R))$ equipped with its Borel sets.
\begin{Borodinstheorem*}[(\cite{Bo1,Bo2})]
There is a probability $\mathbf{P}\in\mathcal P(\bolds{\Om},\bolds
{\mathcal{F}})$ such that
\begin{eqnarray*}
&&\mathbf{P}\bigl(A\x\bigl(D([0,1])\x C_0([0,1]\x\mathbb R)\bigr)\bigr)\\
&&\qquad=\mu(A)\qquad\forall
A\in\B
(\Om);
\\
&&
\mathbf{P}\bigl(\Om\x[(B_\alpha, L_\alpha)\in B]\bigr)\\
&&\qquad=\mathbf{Q}\bigl([(B_\alpha, L_\alpha
)\in B]\bigr)\qquad
\forall B\in\B\bigl(D([0,1])\x C_0([0,1]\x\mathbb R)\bigr);
\end{eqnarray*}
where $\mathbf{Q}=\rmtexttt{dist} (B_\alpha, L_\alpha)\in\mathcal
P(D([0,1])\x
C_0([0,1]\x\mathbb R)$ and such that
\[
(B_{\xi,n}, L_{\xi,n})
\mathop{\lra}^{m}_{n\to\infty} (B_\alpha, L_\alpha) \qquad\text{in }
D([0,1])\x C_0([0,1]\x\mathbb R).
\]
\end{Borodinstheorem*}

This ``coupling version'' is given in~\cite{Bo1} and~\cite{Bo2}. Equivalence with
the distributional version above follows from a general theorem of
Skorokhod; see~\cite{Bi}.

We'll use the following proposition.\vadjust{\goodbreak}
\begin{DCP*}
If $M$ is a metric space, and $\Psi
\dvtx D([0,1])\x C_0([0,1]\x\mathbb R)\to M$ is continuous, then
\[
\Psi(B_{\xi,n}, L_{\xi,n})
\mathop{\lra}^{m}_{n\to\infty} \Psi(B_\alpha, L_\alpha)
\qquad\text{in }
M.
\]
\end{DCP*}
\begin{Localtimelemma*}
For $E\subset\mathbb R$ a finite union of closed, bounded intervals,
\begin{eqnarray*}Y_{E,n}:=\frac1{\overline a_{\xi}(n)}\min_{k\in a_{\xi
}(n)E}N_{n,k}
\mathop{\lra}^{m}_{n\to\infty} \min_{x\in E} L_\alpha(1,x).
\end{eqnarray*}
\end{Localtimelemma*}
\begin{pf}
Since $E\subset\mathbb R$ a finite union of closed, bounded intervals, we
have (using tightness of $\{L_{\xi,n}\dvtx n\in\mathbb N\}$ in $C_0([0,1]\x
\mathbb R)$) that
\[
Y_{E,n}-\min_{x\in E} L_{\xi,n}(1,x)
\mathop{\lra}^{m}_{n\to\infty} 0.
\]

The function $F\dvtx D([0,1])\x C_0([0,1]\x\mathbb R)\to\mathbb R$ defined by
\[
F(X,Y):=\min_{t\in E}Y(t)
\]
is continuous. By the proposition, and
Borodin's theorem,
\begin{eqnarray*}F(B_{\xi,n}, L_{\xi,n})
\mathop{\lra}^{m}_{n\to\infty} F(B_\alpha, L_\alpha)=\min_{x\in E} L_\alpha(1,x).
\end{eqnarray*}
The lemma follows from this.
\end{pf}

\subsection*{Hyperspace}

Let $\mathcal H$ be the hyperspace of all nonempty closed, bounded
subsets of
$\mathbb R$. Equip $\mathcal H$ with the \textit{Hausdorff metric},
\[
h(A, A') : =
\inf\{r > 0 \dvtx A\subset\mathcal N(A',r)\text{ and }A'
\subset\mathcal N(A,r)\}
\]
for $A, B\in\mathcal H$ where, for $A\in\mathcal
H, x\in\mathbb R$ and $r>0$,
\[
\mathcal N(A,r): = \Bigl\{x\in\mathbb R\dvtx \inf_{y\in
A} |x-y|\le r\Bigr\}.
\]
As is well known, $(\mathcal H,h)$ is a locally compact, separable metric
space.

The \textit{range} of the $\mathbb Z$-random walk is $V_n:=\{S_j\dvtx 0\le j\le
n-1\}$. Note that
$\frac1{a(n)}V_n=B_{\xi,n}([0,1])\in\mathcal H$.
\begin{Hyperspaceconvergencelemma*}
Suppose that $(S_1,S_2,\ldots)$ is an aperiodic, $\alpha$-stable random
walk on $\mathbb Z$ with $1<\alpha\le2$, then
\[
\frac1{a(n)}V_n\stackrel{m}\lra\overline{B_{\alpha}([0,1])}
\qquad\text{in }
\mathcal H.
\]
\end{Hyperspaceconvergencelemma*}
\begin{pf}
The function $X\mapsto\overline{X([0,1])}$ is continuous $D([0,1])\to
\mathcal
H$.
\end{pf}

\subsection*{Dyadic partitions and sets}

For $\kappa\in\mathbb N$, let $\D_\kappa$ be the \textit{dyadic partition
of order $\kappa$} defined by
\[
\D_\kappa:=\biggl\{\biggl[\frac{p}{2^\kappa},\frac{p+1}{2^\kappa}\biggr], p\in
\mathbb
Z\biggr\}.
\]

A \textit{closed dyadic set} is a finite union of elements of $\bigcup
_{\kappa\ge1}\D_\kappa$. An \textit{open dyadic set} is the interior of
a closed dyadic set. The \textit{order} of a dyadic set is the minimal
$\kappa\in\mathbb N$ so that the (closure of the) dyadic set is a union of
elements of $\D_\kappa$. Let
$\mathcal D^o_\kappa$ and $\overline{\mathcal D}_\kappa$ denote the
collections of open and closed dyadic sets of order $\kappa$, respectively.

For $E\subset\mathbb R$ bounded, nonempty and $\kappa\ge1$ let:

$\bullet$ $C_\kappa(E)$ be the largest closed dyadic set of order $\kappa
$ contained in $E^o$;

$\bullet$ $U_\kappa(E)$ be the smallest open dyadic set of order $\kappa$
containing $E$. Note that $C_\kappa(E)\subset U_\kappa(E)\ne
\varnothing$ and that it is possible that $C_\kappa(E)=\varnothing$.

For $\kappa\in\mathbb N, \Upsilon\in\mathcal D^o_\kappa, \G\in
\overline
{\mathcal D}_\kappa$ satisfying $\G\subset\Upsilon$, define the set
\[
\mathcal U(\kappa,\G,\Upsilon):=\{E\in\mathcal H\dvtx C_\kappa(E)=\G
,U_\kappa(E)=\Upsilon\}.
\]
These sets are not open in $\mathcal H$, but are Borel sets in
$\mathcal H$ with the additional property that

\Smi\mbox{ }$\forall E\in\mathcal U(\kappa,\G,\Upsilon) \exists
\delta>0$ such that $
\mathcal N(E,\delta)\subset\Upsilon$.

The sets
$\mathcal U(\kappa,\G,\Upsilon)$ and $\mathcal U(\kappa,\G
',\Upsilon')$
are disjoint unless $\G=\G'$ and $\Upsilon=\Upsilon'$.

\subsection*{Admissibility}

For $\mathcal E>0$, we call a pair $(\G,\Upsilon
)\in
\bigcup_{\kappa\ge1}\overline{\mathcal D}_\kappa\x\mathcal D^o_\kappa$
\textit{$\mathcal E$-admissible} if:

\begin{longlist}
\item
$\mu=\mu(\G,\Upsilon):=\rmtexttt{Leb} (\Upsilon
\setminus\G
)<\mathcal E$;

\noindent and for $N\in\mathbb N$, we call $(\G,\Upsilon)$ $(N,\mathcal E)$-\textit{admissible}
if in addition

\item
$MH(3\mu)+3\mu\log N<\mathcal E$ where $M=M(\G,\Upsilon
):=\rmtexttt{Leb} (\Upsilon)$ and $H(t):=-t\log t-(1-t)\log(1-t)$.
\end{longlist}

Note that if $A\subset\mathbb R$ is Riemann integrable, then $\forall
\mathcal E>0, \exists$ a
$(N,\mathcal E)$-admissible pair $(\G,\Upsilon)\in\overline{\mathcal
D}_\kappa
\x\mathcal D^o_\kappa$ so that
$\G\subset A \subset\Upsilon$ and that in this case $(C_\kappa
(A),U_\kappa(A))$ is also $(N,\mathcal E)$-admissible.
\begin{lemma}\label{lemma4}
For each $\mathcal E>0, N\in\mathbb N, \e>0, \exists \kappa\in
\mathbb N$
and $\theta>0$ and a finite collection of $(N,\mathcal E)$-admissible
pairs
\[
\{(\G_j,\Upsilon_j)\}_{j\in J}\subset\overline{\mathcal D}_\kappa\x
\mathcal
D^o_\kappa
\]
satisfying:
\begin{longlist}
\item
\[
\mathbf{P} \biggl(\biguplus_{j\in J}G_j\biggr)>1-\e
\qquad\text{where } G_j:=\bigl[\overline{B_\alpha([0,1])}\in\mathcal U(\kappa,\G
_j,\Upsilon_j)\bigr].
\]
For large enough $n\ge1$,%\vspace*{1pt}

\item
\[
\mathbf{P} (G_{j,\theta,n})> (1-\varepsilon
)\cdot
\mathbf{P} (G_j)\qquad\forall j\in J,\vadjust{\goodbreak}
\]
where
\[
G_{j,\theta
,n}:=\biggl[\frac1{\overline{a}(n)}\min_{k\in a(n)\G_j}N_{n,k}>\theta,
\frac
1{a(n)}V_n\subset\Upsilon_j\biggr]\cap G_j.
\]
\end{longlist}
\end{lemma}
\begin{pf}
By Lemma~\ref{lemma3}, $B_\alpha([0,1])$ is a.s. Riemann
integrable, so $\exists \kappa\in\mathbb N$ and a finite collection of
$(N,\mathcal E)$-admissible pairs $\{(\G_j,\Upsilon_j)\}_{j\in
J}\subset
\overline{\mathcal D}_\kappa\x\mathcal D^o_\kappa$ satisfying~(i).

Suppose that $\mathbf{P} (G_j)\ge\eta>0\mbox{  }\forall j\in J$.

By Lemma~\ref{lemma3}, $\min_{x\in\G_j} L_\alpha(1,x)>0$ a.s. on $G_j$
$\forall j\in
J$. This and (\smiley) ensure that $\exists \theta>0$ such that
$\forall
j\in J$,
\begin{eqnarray*}\mathbf{P} \Bigl(\Bigl[\min_{x\in\G_j} L_\alpha(1,x)>2\theta
\Bigr]\cap
\bigl[\mathcal N(\overline{B_\alpha([0,1])},\theta)\subset\Upsilon_j\bigr]\cap G_j\Bigr)>
\biggl(1-\frac
{\varepsilon}2\biggr)
\mathbf{P} (G_j).
\end{eqnarray*}
By the local time, and hyperspace convergence lemmas, for $n\ge1$ large
\begin{eqnarray*}
\mathbf{P} \biggl(\biggl[h\biggl(\frac1{a(n)}V_n,\overline{B_\alpha([0,1])}\biggr)\ge\theta\biggr]\biggr)&<&\frac
{\eta\e}4;
\\
\mathbf{P} \biggl(\biggl[\biggl|\min_{x\in\G_j} L_\alpha(1,x)-\frac1{a(n)}\min_{k\in
a_{\xi
}(n)\G_j}N_{n,k}\biggr|>\theta\biggr]\biggr)&<&\frac{\eta\e}4.
\end{eqnarray*}
Statement (ii) follows from this.
\end{pf}

%s4 ###
\section{Relative complexity of RWRS}\label{sec4}

We prove Theorem~\ref{pi3}(1).

Fix a finite, $S$-generator $\beta\in\frak P(Y,\C,\mu)$, and let
$P=P_\beta\in
\frak P(Z,\B,m)$ defined by
$P(x,y):=\alpha(x)\x\beta(y)$ where $\alpha(x):=[x_0]$, then
\[
P_0^{n-1}(T)(x,y)= \alpha_0^{n-1}(R)(x)\x\beta
_{V_n(x)}(S)(y),
\]
where
\[
\beta_{V_n(x)}(S):=\bigvee
_{k\in
V_n(x)}S^{-k}\beta.
\]

Define for $n\in\mathbb N, \e>0$ [as in the proof of~(\ref{bstar})],
$\Pi
_n\dvtx\Om\to2^{P_0^{n-1}(T)}$ by
\[
\Pi_n(x):=\bigl\{a\in P_0^{n-1}(T)\dvtx m\bigl(a\|\B(\Om)\x Y\bigr)(x)>0\bigr\}.
\]
Note that for fixed $x\in\Om$, if $z\in\Pi_n(x)$, then $z$ is of
form
\[
z=(x_0^{n-1},w):=[x_0^{n-1}]\x\bigvee_{j\in V_n(x)}S^{-j}w_j\mbox{
}(w_j\in\beta).
\]

Now define $\Phi_{n,\e}, \mathcal Q(P,n,\e)\dvtx\Om\to\mathbb N$ by
\begin{eqnarray*}
\Phi_{n,\e}(x)&:=&\min \biggl\{\# F\dvtx F\subset\Pi_n(x)\dvtx m\biggl(\bigcup_{a\in
F}a\|
\B(\Om)\x Y\biggr)(x)>1-\e\biggr\};
\\
\mathcal Q(P,n,\e)(x)&:=&\max \bigl\{\#\{c\in\Pi_n(x)\dvtx \overline d_n(a,c)\le\e
\}\dvtx a\in\Pi_n(x)\bigr\},
\end{eqnarray*}
where $\overline d_n$ is the $(T,P,n)$-Hamming metric
\[
\overline d_n([a_0,\ldots,a_{n-1}],[c_0,\ldots,c_{n-1}])=\frac1n\#\{0\le
k\le
n-1\dvtx a_k\ne c_k\}.
\]
As before,
\renewcommand{\theequation}{\Bicycle}
\begin{equation}
\frac{\Phi_{n,\e}(x)}{\mathcal
Q(P,n,\e
)(x)}\le K^{(T)}_{\B(\Om)\x Y}(P,n,\e)(x)\le\Phi_{n,\e}(x).
\end{equation}

To establish Theorem~\ref{pi3}(1),
it suffices by~(\ref{bicycle}) to show that
\renewcommand{\theequation}{\Football}
\begin{equation}\frac1{a(n)}\log_2\Phi_{n,\e}
\mathop{\lra}^{{\frak d}}_{n\to\infty, \e\to 0}\rmtexttt{Leb}
(B_\alpha
([0,1])) h(S,\beta)
\end{equation}
and
\renewcommand{\theequation}{\protect
\includegraphics{688i04.eps}
}
\begin{equation}
\frac1{a(n)}\log_2\mathcal Q(P,n,\e)
\mathop{\lra}^{m}_{n\to\infty, \e\to0} 0.
\end{equation}\vspace*{6pt}

\textit{Proof of}~(\ref{football}).\quad
In order to use Lemma~\ref{lemma4}, we consider $(\mathbf{Z},\bolds{\B}(\mathbf{Z}),
\mathbf{m},\mathbf{T})$ where
\begin{eqnarray*}
\mathbf{Z}&:=&\bolds{\Om}\x Y\cong Z\x D([0,1])\x
C_0([0,1]\x\mathbb R),\\
\mathbf{m}&:=&\mathbf{P}\x\nu,\\
\mathbf{T}(x,y,t)&:=&(T(x,y),t)
\end{eqnarray*}
and prove that on $(\mathbf{Z},\bolds{\B}(\mathbf{Z}),\mathbf{m})$,
\renewcommand{\theequation}{$\hat{\text{\Football}}$}
\begin{equation}\label{wdfootball}
\frac1{a(n)}\log_2\Phi_{n,\e}
\mathop{\lra}^{m}_{n\to\infty, \e\to 0} \rmtexttt{Leb} (B_\alpha([0,1]))
h(S,\beta);
\end{equation}
[which implies~(\ref{football}) on $({Z,\B(Z),m})$], deducing
(\ref{wdfootball}) from
\renewcommand{\theequation}{$\tilde{\text{\Football}}$}
\begin{equation}\label{wdtldfootball}
\frac1{a(n)}I\bigl(P_0^{n-1}(T)\|\B(\Om)\x Y\bigr)
\mathop{\lra}^{m}_{n\to\infty} \rmtexttt{Leb} (B_\alpha([0,1]))
h(S,\beta),
\end{equation}
where $I(\alpha\|\C)$ is \textit{conditional information} defined by
\[
I(\alpha\|\C)(x):=\log\frac1{\mathbf{m}(\alpha(x)\|\C)(x)}.
\]
\textit{Proof of}~(\ref{wdtldfootball}).\quad
By the above
\[
P_0^{n-1}(T)(x,y)= \alpha_0^{n-1}(R)(x)\x\beta_{V_n(x)}(S)(y),
\]
whence
\[
I\bigl(P_0^{n-1}(T)\|\B(\Om)\x Y\bigr)(x)=\log\frac1{\nu(\beta
_{V_n(x)}(S)(y))}=I\bigl(\beta_{V_n(x)}(S)\bigr)(y).
\]
The idea of the proof is to approximate $V_n(x)$ with sequences of sets
of form $F_{\Lambda,n}:=(a(n)\Lambda)\cap\mathbb Z$ where $\Lambda
\subset
\mathbb R$ is a finite union of disjoint, bounded, intervals.\vadjust{\goodbreak}

Any such sequence $\{F_{\Lambda,n}\dvtx n\ge1\}$
satisfies \ittexttt{F\o llner's condition}:
\[
\frac{\#(F_{\Lambda,n}\D(F_{\Lambda,n}+j))}{\# F_{\Lambda,n}}
\mathop{\lra}_{n\to\infty} 0\qquad\forall j\in\mathbb Z.
\]

Moreover,
\[
\# F_{\Lambda,n}=(a(n)\Lambda)\cap\mathbb Z=a(n)\rmtexttt{Leb}
(\Lambda
)\pm2M,
\]
where $\Lambda$ is a union of $M$ disjoint intervals.

Thus, by Kieffer's Shannon--MacMillan theorem (\cite{Ki}---see also~\cite{MO})
\renewcommand{\theequation}{\mbox{SM}}
\begin{equation}\label{SM}
\frac1{a(n)}I(\beta_{F_{\Lambda,n}}(S))
\mathop{\lra}^{m}_{n\to\infty} h(S,\beta)\rmtexttt{Leb} (\Lambda).
\end{equation}

Fix $\e>0$ and
let
$\kappa\in\mathbb N, \theta>0$ and the finite collection  $\{(\G
_j,\Upsilon_j)\}_{j\in J}\subset\overline{\mathcal D}_\kappa\x\mathcal
D^o_\kappa$
%be as advertised by
of $\e$-admissible pairs be as in
Lemma~\ref{lemma4}.

By~(\ref{SM}) for large $n\ge1, \exists H_n\in\B(Y)$ so that
$\nu(H_n)>1-\e$ and such that $\forall y\in H_n, \Lambda\in\{\G
_j,\Upsilon
_j\}_{j\in J}$,
\[
\frac1{a(n)}I(\beta_{F_{\Lambda,n}}(S))(y)=(1\pm\e)h(S,\beta
)\rmtexttt{Leb}
(\Lambda).
\]

For $\omega\in G_{j,\theta,n}$,
\[
\frac1{a(n)}V_n(\omega)\subset\mathcal N(\overline{B_\alpha([0,1])},\theta
)\subset
\Upsilon_j
\]
and
\[
(a(n)\G_j)\cap\mathbb Z\subset\{k\in\mathbb Z\dvtx
N_{n,k}(\omega)>\theta\overline{a}(n)\}\subset V_n(x).
\]
Thus
\[
F_{\G_j,n}\subset V_n(\omega)\subset F_{\Upsilon_j,n},
\]
whence for $y\in H_n$,
\begin{eqnarray*}
&&
h(S,\beta)\rmtexttt{Leb} (B_\alpha([0,1]))-\e\\
&&\qquad<\frac
1{a(n)}I(\beta_{\Lambda_{n,\G_j}}(S))
\le\frac1{a(n)}I(\beta
_{ V_n(\omega
)}(S))\\
&&\qquad\le\frac1{a(n)}I(\beta_{\Lambda_{n,\Upsilon_j}}(S))\\
&&\qquad<h(S,\beta
)\rmtexttt{Leb} (B_\alpha([0,1]))+\e.
\end{eqnarray*}

Thus
\renewcommand{\theequation}{$\tilde{\text{\Football}}$}
\begin{eqnarray}
&&\mathbf{P}\biggl(\biggl[\frac1{a(n)}
I\bigl(\beta_{V_n(\omega)}(S)\bigr)=h(S,\beta)\rmtexttt{Leb} (B_\alpha
([0,1]))\pm\e\biggr]\biggr)\nonumber\\
&&\qquad\ge
\sum_{j\in J}\mathbf{P}(G_{j,\theta,n})\nu(H_n)\nonumber\\[-8pt]\\[-8pt]
&&\qquad>(1-\e)^2\sum
_{j\in
J}\mathbf{P}(G_{j})\nonumber\\
&&\qquad>(1-\e)^3.\nonumber
\end{eqnarray}
\begin{Remark*} We note that the methods of the proof of
(\ref{wdtldfootball}) can be adapted to prove Theorem 7.1
in~\cite{LeG-R}, namely
\[
\frac{\# V_n}{a(n)} \stackrel{\frak d}\lra \rmtexttt{Leb} (B_\alpha
([0,1])).
\]
\end{Remark*}
\textit{Proof of}~(\ref{wdfootball}).\quad
By~(\ref{wdtldfootball}), $\forall \e>0, \exists
N_\e$ such that $\forall n>N_\e \exists
G_n\in\B(\bolds{\Om})$ so that for $x\in G_n$,
\[
\nu(H_{n,x})>1-\e,
\]
where
\[
H_{n,x}:=\bigl\{y\in Y\dvtx
\nu(\beta_{V_n(x)}(S)(y))=e^{-a(n)\rmtexttt{Leb}
(B_\alpha([0,1])h(S,\beta)(1\pm\e))}\bigr\}.
\]
Let $F_{n,\e,x}:=\{\beta_{V_n(x)}(S)(y)\dvtx y\in H_{n,x}\}$.
It follows that
\[
\log\# F_{n,\e,x}=a(n)\rmtexttt{Leb} \bigl(B_\alpha([0,1])h(S,\beta)(1\pm
\e)\bigr).
\]
Thus
\[
\log\Phi_{n,\e}(x)\le a(n)\rmtexttt{Leb} \bigl(B_\alpha([0,1])h(S,\beta
)(1+\e)\bigr).
\]
On the other hand, if
$F\subset\Pi_n(x), m(\bigcup_{a\in F}a\|\B(\Om)\x Y)(x)>1-\e$, then
$F\supset F_{n,2\e,x}$, whence
\renewcommand{\theequation}{$\hat{\text{\Football}}$}
\begin{equation}
\log\Phi_{n,\e}(x)\ge a(n)\rmtexttt{Leb} \bigl(B_\alpha([0,1])h(S,\beta
)(1-2\e)\bigr).
\end{equation}

\textit{Proof of} \mbox{(\ref{radioactivity})}.\quad
Fix $\e=\mathcal E>0$. Let $\kappa\in\mathbb N$ and
$\theta>0$
and the finite collection of $(\#\beta,\mathcal E)$-admissible pairs
\[
\{(\G_j,\Upsilon_j)\}_{j\in J}\subset\overline{\mathcal D}_\kappa\x
\mathcal
D^o_\kappa
\]
be as in Lemma~\ref{lemma4}.

(1) For large $n, x\in G_{j,n,\theta}, a=(x_0^{n-1},w),
a'=(x_0^{n-1},w')\in\Pi_n(x)$,
\[
\#\{i\in V_n(x)\dvtx w_i\ne w_i'\}\le
a(n)\biggl({\mu_j}+\frac{d_n(a,a')}\theta\biggr),
\]
where
$\mu_j:=\mu(\G_j,\Upsilon_j)=\rmtexttt{Leb} (\Upsilon_j\setminus
\G_j)$.
\begin{pf}
Let $x\in G_{j,n,\theta}$, and let
\[
K_n(x):=\{i\in V_n(x)\dvtx w_i\ne w_i'\}.
\]
Then since
\[
\Lambda_{n,\G_j}\subset V_n(x)\subset\Lambda_{n,\Upsilon_j}
\]
we have that
\[
K_n\subset K_n\cap\Lambda_{n,\G_j}\cup\Lambda_{n,\Upsilon
_j}\setminus
\Lambda_{n,\G_j},
\]
whence, for large $n$,
\[
\# K_n\le\#(K_n\cap\Lambda_{n,\G_j})+\#(\Lambda_{n,\Upsilon
_j}\setminus
\Lambda_{n,\G_j})
\le\#(K_n\cap\Lambda_{n,\G_j})+\mu_j a(n).
\]

Now,
\begin{eqnarray*} \#(K_n\cap\Lambda_{n,\G_j})& \le&
\frac1{\theta\overline a(n)}\sum_{k\in\Lambda_{n,\G_j}}N_{n,k}1_{K_n}(k)\\
&=&
\frac1{\theta\overline a(n)}\sum_{k=0}^{n-1}\#\bigl\{0\le i\le n-1\dvtx w_{s_i(x)}\ne
w'_{s_i(x)}\bigr\}\\
&=&\frac{n}{\theta\overline a(n)}  d_n(a,a')\\
&\lesssim&
\frac
1{\theta} a(n) d_n(a,a').
\end{eqnarray*}
\upqed\end{pf}

(2) For $n$ large,
\[
\max_{x\in G_{j,n,\theta}}\mathcal Q\biggl(P,n,\frac{\mu_j}\theta\biggr)(x)\le
e^{\mathcal Ea(n)(1+o(1))}.
\]
\begin{pf}
Fix $x\in G_{j,n,\theta}, z=(x_0^{n-1},u)\in\Pi_n(x)$, then
\begin{eqnarray*}
&&\biggl\{a\in\Pi_n(x)\dvtx a\subset B\biggl(n,P,z,\frac{\mu
_j}\theta
\biggr)\biggr\}\\
&&\qquad
\mathop{\subseteqq}_{\sim}\biggl\{v\in\beta^{V_n}\dvtx
d_n((x_0^{n-1},u),(x_0^{n-1},v))<\frac{\mu_j}\theta\biggr\}\\
&&\qquad\stackrel{(1)}{\subseteq}\bigl\{v\in\beta^{V_n}\dvtx \#\{i\in V_n(x)\dvtx
v_i\ne
u_i\}\le2\mu_j a(n)\bigr\}.
\end{eqnarray*}
Thus for $n$ large,
\begin{eqnarray*}
&&\# \biggl\{\Pi_n(x)\dvtx a\subset B\biggl(n,P,z,\frac{\mu_j}\theta
\biggr)\biggr\}
\le
\pmatrix{\# V_n(x)\cr2\mu_j a(n)}|\beta|^{2\mu_j a(n)}\\
&&\qquad\le
\pmatrix
{M_ja(n)\cr2\mu_j a(n)}|\beta|^{2\mu_j a(n)} \qquad\text{where } M_j:=\rmtexttt{Leb} (\Upsilon_j);\\
&&\qquad\le
e^{M_jH(2\mu_j)a(n)(1+o(1))}|\beta|^{2\mu_j a(n)}
\qquad\text{by Stirling's formula};\\
&&\qquad=e^{(M_jH(2\mu_j)+2\mu_j\log|\beta
|)a(n)(1+o(1))}\\
&&\qquad= e^{\mathcal Ea(n)(1+o(1))}.
\end{eqnarray*}
\upqed\end{pf}

By (2), if $\delta=\delta(\mathcal E):=\min_{j\in J}\frac{\mu
_j}\theta$, then $\delta
>0$ and
\renewcommand{\theequation}{\protect
\includegraphics{688i04.eps}
}
\begin{equation}
\mathbf{P}\bigl([\log_2\mathcal Q(P,n,\delta)<\mathcal E a(n)]\bigr)>\sum
_{j\in
J}\mathbf{P}(G_{j,\theta,n})>(1-\e)^2.
\end{equation}
As mentioned above, this establishes Theorem~\ref{pi3}.

\section*{Concluding remarks and questions}

Recently in~\cite{BK}, Borodin's theorem~\cite{Bo1} (coupling version) has been
established for strongly aperiodic random walks driven by Markov chains, and Theorem 2 can now
be proven with the same methods in this case.

However, Theorem 2 applies neither to a RWRS whose jump random variables are $1$-stable
nor to a generalized RWRS over $\mathbb Z^2$ whose jump random variables are
centered and in the domain of attraction of standard normal
distribution on $\mathbb R^2$. Other methods are needed to treat these
cases due to the lack of ``smooth local time'' of the relevant limit processes.

It is still conceivable that in both cases there are $1$-regularly
varying relative complexity sequence whence (or otherwise)
\[
{\rmtexttt{E\mbox{-}dim}} (T,\text{Base})=1.
\]

Nothing is known about the relative complexity of generalized RWRSs
over continuous groups (as in~\cite{Ba}) or of ``smooth RWRSs'' (as in~\cite{Ru}).

\section*{Acknowledgments}
The author would like to acknowledge helpful conversations with Omer
Adelman, Nathalie Eisenbaum and Haya Kaspi. Lemma~\ref{lemma3} is due to Haya
Kaspi and Nathalie Eisenbaum and reproduced here with their kind
permission.

%suskaldyti doi

% imsref loaded by lrinkeviciute, 2011-08-16 16:20:43

%
\printaddresses

\end{document}